\documentstyle{amsart}
\font\bbb = msbm10

\def\Bbb#1{\hbox{\bbb #1}}

\def\R{{\Bbb R}}
\def\C{{\Bbb C}}

\def\D{{\Bbb D}}

\newcommand{\dis}{\displaystyle}

\def\DATE{
}

\newtheorem{theorem}{Theorem}
\newtheorem{proposition}{Proposition}[section]
\newtheorem{lemma}[proposition]{Lemma}

\begin{document}

\large

\title[Solvability of the generalized Possio equation]
{Solvability of the generalized Possio equation in 2D subsonic aeroelasticity.}
\author[Peter L. Polyakov]{Peter L. Polyakov}
\address{Department of Mathematics, University of Wyoming, Laramie, WY 82071, USA}
\email{ polyakov@@uwyo.edu}
\date{\today}

\begin{abstract}
We study solvability of the {\it generalized Possio integral equation} - a tool in analysis of
a boundary value problem in 2D subsonic aeroelasticity with the Kutta-Joukowski
condition - {\it "zero pressure discontinuity"} -  $\psi(x,0,t)=0$
on the complement of a finite interval in the whole real line $\R$. The corresponding
problem with boundary condition on finite intervals adjacent to the "chord" was considered
in \cite{P}.
\end{abstract}
\maketitle

\centerline{\DATE}

\bigskip

\section{Introduction.}\label{Introduction}
\indent
We consider the linearized subsonic inviscid compressible flow equation in 2D
(\cite{Se}, \cite{BAH}, \cite{Ba2})
\begin{equation}\label{LinearEquation}
a^2\left(1-M^2\right)\frac{\partial^2\phi}{\partial x^2}
+a^2\frac{\partial^2\phi}{\partial y^2}
=\frac{\partial^2\phi}{\partial t^2}
+2Ma\frac{\partial^2\phi}{\partial t\partial x},
\end{equation}
where $a$ is the speed of sound,
${\dis M=\frac{U}{a}<1}$ - the Mach number, $U$ - free stream velocity,
$\phi(x,y,t)$ - small disturbance velocity potential, considered on
$$\left\{(x,y,t):(x,y)\in \R^2\setminus\left\{x\in [-1,1],\ y=0\right\},\
0\leq t<\infty\right\},$$
with boundary conditions:
\begin{itemize}
\item
flow tangency condition
\begin{equation}\label{FlowTangency}
\frac{\partial\phi}{\partial y}(x,0,t)=w(x,t),\ |x|<b,
\end{equation}
where $b$ is the "half-chord", and $w$ is the given normal velocity of the wing,
without loss of generality we will assume in what follows that $b=1$,
\item
{\it "strong Kutta-Joukowski condition"} for the acceleration potential
$$\psi(x,y,t):=\frac{\partial\phi}{\partial t}+U\frac{\partial\phi}{\partial x},$$
\begin{equation}\label{Kutta}
\psi(x,0,t)=0\ \mbox{for}\ |x|>1,
\end{equation}
\item
far field condition
$$\phi(x,y,t)\to 0,\ \mbox{as}\ |x|\to\infty,\ \mbox{or}\ y\to\infty.$$
\end{itemize}

\indent
In \cite{P} we considered a weaker form of the boundary condition (\ref{Kutta}),
requiring for acceleration potential to be zero not on the whole $\R\setminus[-1,1]$,
but only on finite intervals adjacent to the interval $[-1,1]$. In the present paper we address
the problem with "auxiliary boundary conditions" from (\cite{BAH}, p. 319) on the whole
$\R\setminus[-1,1]$, combining technical tools from \cite{P} with Possio's approach as it is described in
the fascinating book \cite{BAH}. We also generalize the Possio's scheme of construction of
a solution of the boundary value problem (\ref{LinearEquation}), (\ref{FlowTangency}), (\ref{Kutta})
by including the case of arbitrary (not only harmonic) dependence of
$w(x,t)$ (and consequently of a solution) on time.\\
\indent
Before formulating the main result of the article we introduce a couple of notations.
For a function $w(x,t)$ on $[-1,1]\times \R^{+}$ we denote by
$${\widehat w}(x,s)=\int_0^{\infty}e^{-st}w(x,t)dt,$$
the Laplace transform of $w(x,t)$, an analytic function of $s=\sigma+i\nu$ in the half-plane
$\left\{s:\mbox{Re}s>\sigma_0\right\}$.
In section \ref{GOperatorResolvent} we introduce function $D_{N}(s)$ (formula (\ref{D_P})
and Proposition~\ref{GResolvent}), analytic in the same half-plane
and depending only on function $H_0^{(1)}$ - the Hankel function of the first kind of order $0$.\\
\indent
The main result of the article is the theorem below.

\begin{theorem}\label{Main}\ Let function $D_{N}(s)$ from formula (\ref{D_P}),
mentioned above, have no zeros in the strip $\left\{s:\mbox{Re}s\in [\sigma_1,\sigma_2]\right\}$,
where $\sigma_1>\sigma_0$. Let $w(\cdot,t) \in L^2[-1,1]$ be such that for some $\epsilon>0$
\begin{equation}\label{wCondition}
{\dis \left\|\widehat{w}(\cdot,\sigma+i\nu)\right\|_{L^2[-1,1]}
<\exp\left\{-e^{|\nu|}\cdot(1+|\nu|)^{4+\epsilon}\right\}\ \mbox{for}\
\sigma\in[\sigma_1,\sigma_2] }
\end{equation}
\indent
Then equation (\ref{LinearEquation}) has a solution of the form
\begin{equation}\label{Solution}
\phi(x,y,t)=-\frac{1}{2\pi i}
\int_{\sigma^{\prime}-i\infty}^{\sigma^{\prime}+i\infty}dse^{s(t+cx)}\frac{e^{\lambda(s)x}}{U}
\int_{-1}^1p(\xi,s)d\xi\int_{-\infty}^x e^{-\lambda(s)u}
\left(\frac{\partial H_0^{(1)}}{\partial y}(\zeta)\Bigg|_{\eta=0}\right)du,
\end{equation}
where ${\dis c=\frac{M}{a(1-M^2)} }$, ${\dis \lambda(s)=-\frac{s\left(1+cU\right)}{U} }$,
$H_0^{(1)}$ - the Hankel function of the first kind of order $0$, and
$$\zeta=\frac{is}{a\sqrt{1-M^2}}\sqrt{\frac{(u-\xi)^2}{1-M^2}+(y-\eta)^2}.$$
This solution is independent of $\sigma^{\prime}\in [\sigma_1,\sigma_2]$,
satisfies boundary conditions above, and function $p(\xi,s)$ satisfies the estimate
\begin{equation}\label{pEstimate}
\int_{-1}^1\left|p(\xi,\sigma^{\prime}+i\nu)\right|^pd\xi<C
\end{equation}
with $C>0$ independent of $s$ for arbitrary ${\dis p<\frac{4}{3} }$.
\end{theorem}

\section{Particular solution of equation (\ref{LinearEquation}).}\label{Particular}

\indent
In this section we follow Possio's idea of constructing a special solution
of equation (\ref{LinearEquation}) that represents the acceleration potential. It is easy
to see that since equation (\ref{LinearEquation}) is linear, the acceleration potential, defined
in the linear model of subsonic flow by formula in (\ref{Kutta}) also satisfies equation
(\ref{LinearEquation}). In order to construct this special solution
of differential equation (\ref{LinearEquation}) we use two lemmas below.

\begin{lemma}\label{ReducedWave}\ If function $\Psi(x,y)$ satisfies the
reduced wave equation
\begin{equation}\label{RedWave}
a^2\left(1-M^2\right)\frac{\partial^2\Psi}{\partial x^2}(x,y)
+a^2\frac{\partial^2\Psi}{\partial y^2}(x,y)-\frac{s^2}{1-M^2}\Psi(x,y)=0,
\end{equation}
then
function
\begin{equation}\label{psi_Psi}
\psi(x,y,t)=\Psi(x,y)e^{s(t+cx)}
\end{equation}
with
\begin{equation}\label{cFormula}
c=\frac{M}{a(1-M^2)}
\end{equation}
satisfies equation (\ref{LinearEquation}).
\end{lemma}
\indent
{\bf Proof.}\ For $\psi$ defined in (\ref{psi_Psi}) we have
$$a^2\left(1-M^2\right)\frac{\partial^2\psi}{\partial x^2}
+a^2\frac{\partial^2\psi}{\partial y^2}
-\frac{\partial^2\psi}{\partial t^2}
-2Ma\frac{\partial^2\psi}{\partial t\partial x}$$
$$=e^{s(t+cx)}\left[a^2\left(1-M^2\right)
\left(\frac{\partial^2\Psi}{\partial x^2}(x,y)
+2sc\frac{\partial\Psi}{\partial x}(x,y)+s^2c^2\Psi(x,y)\right)
+a^2\frac{\partial^2\Psi}{\partial y^2}(x,y)\right.$$
$$\left.-s^2\Psi(x,y)-2Ma\left(s\frac{\partial\Psi}{\partial x}(x,y)
+s^2c\Psi(x,y)\right)\right]$$
$$=e^{s(t+cx)}\left[a^2\left(1-M^2\right)\frac{\partial^2\Psi}{\partial x^2}(x,y)
+a^2\frac{\partial^2\Psi}{\partial y^2}(x,y)\right.$$
$$\left.+\left(a^2\left(1-M^2\right)s^2c^2-2Mas^2c-s^2\right)\Psi(x,y)\right]$$
$$=e^{s(t+cx)}\left[a^2\left(1-M^2\right)\frac{\partial^2\Psi}{\partial x^2}(x,y)
+a^2\frac{\partial^2\Psi}{\partial y^2}(x,y)-\frac{s^2}{1-M^2}\Psi(x,y)\right]=0.$$
\qed

\begin{lemma}\label{HankelSolution}\ Let
$$H_0^{(1)}=J_0+iY_0$$
be the Hankel function of the first kind of order $0$.\\
\indent
Then for fixed $\xi$ and $\eta$ function
$$\Psi_{\xi,\eta}(x,y)=\frac{\partial}{\partial \eta}H_0^{(1)}\left(z\right)$$
with
\begin{equation}\label{zFormula}
z=\frac{is}{a\sqrt{1-M^2}}\sqrt{\frac{(x-\xi)^2}{1-M^2}+(y-\eta)^2}
\end{equation}
satisfies equation (\ref{RedWave}).
\end{lemma}
\indent
{\bf Proof.}\ Because of linearity of equation (\ref{RedWave}) it suffices to prove
that function $H_0^{(1)}(z)$ satisfies this equation. For this function we have
$$a^2\left(1-M^2\right)\frac{\partial^2 H_0^{(1)}(z)}{\partial x^2}
+a^2\frac{\partial^2 H_0^{(1)}(z)}{\partial y^2}-\frac{s^2}{1-M^2}H_0^{(1)}(z)$$
$$=a^2\left(1-M^2\right)\frac{\partial}{\partial x}
\left[\frac{dH_0^{(1)}}{dz}(z)\frac{\partial z}{\partial x}\right]
+a^2\frac{\partial}{\partial y}
\left[\frac{dH_0^{(1)}}{dz}(z)\frac{\partial z}{\partial y}\right]
-\frac{s^2}{1-M^2}H_0^{(1)}(z)$$
$$=-a^2\left(1-M^2\right)\frac{s^2}{a^2(1-M^2)}\frac{d^2H_0^{(1)}}{dz^2}(z)
\cdot\left[\left(\sqrt{\frac{(x-\xi)^2}{1-M^2}+(y-\eta)^2}\right)^{-1}
\cdot\frac{x-\xi}{1-M^2}\right]^2$$
$$+\frac{isa^2\left(1-M^2\right)}{a\sqrt{1-M^2}}\frac{dH_0^{(1)}}{dz}(z)\cdot
\left[\frac{1}{1-M^2}\left(\sqrt{\frac{(x-\xi)^2}{1-M^2}+(y-\eta)^2}\right)^{-1}\right.$$
$$\left.-\left(\frac{x-\xi}{1-M^2}\right)^2\left(\sqrt{\frac{(x-\xi)^2}{1-M^2}+(y-\eta)^2}\right)^{-3}
\right]$$
$$-a^2\frac{s^2}{a^2(1-M^2)}\frac{d^2H_0^{(1)}}{dz^2}(z)
\cdot\left[\left(\sqrt{\frac{(x-\xi)^2}{1-M^2}+(y-\eta)^2}\right)^{-1}
\cdot(y-\eta)\right]^2$$
$$+a^2\frac{is}{a\sqrt{1-M^2}}\frac{dH_0^{(1)}}{dz}(z)\cdot
\left[\left(\sqrt{\frac{(x-\xi)^2}{1-M^2}+(y-\eta)^2}\right)^{-1}\right.$$
$$\left.-\left(\sqrt{\frac{(x-\xi)^2}{1-M^2}+(y-\eta)^2}\right)^{-3}\cdot(y-\eta)^2\right]$$
$$-\frac{s^2}{1-M^2}H_0^{(1)}(z)$$
$$=-\frac{s^2}{1-M^2}\frac{d^2H_0^{(1)}}{dz^2}(z)
\left[\left(\frac{(x-\xi)^2}{1-M^2}+(y-\eta)^2\right)^{-1}
\cdot\left(\frac{(x-\xi)^2}{1-M^2}+(y-\eta)^2\right)\right]$$
$$+2\frac{isa}{\sqrt{1-M^2}}\frac{dH_0^{(1)}}{dz}(z)\cdot
\left(\sqrt{\frac{(x-\xi)^2}{1-M^2}+(y-\eta)^2}\right)^{-1}$$
$$-\frac{isa}{\sqrt{1-M^2}}\frac{dH_0^{(1)}}{dz}(z)\cdot
\left(\sqrt{\frac{(x-\xi)^2}{1-M^2}+(y-\eta)^2}\right)^{-1}
-\frac{s^2}{1-M^2}H_0^{(1)}(z)$$
$$=\frac{-s^2}{1-M^2}\left[\frac{d^2H_0^{(1)}}{dz^2}(z)
+\frac{dH_0^{(1)}}{dz}(z)\left(\frac{is}{a\sqrt{1-M^2}}
\sqrt{\frac{(x-\xi)^2}{1-M^2}+(y-\eta)^2}\right)^{-1}+H_0^{(1)}(z)\right]=0,$$
where in the last equality we used the fact that function $H_0^{(1)}$
satisfies the Bessel equation of order zero.\\
\qed

\indent
From Lemma~\ref{HankelSolution} we obtain that for arbitrary $\xi\in[-1,1]$ function
\begin{equation}\label{DoubletPsi}
\Psi_{\xi,0}(x,y)=\frac{\partial}{\partial \eta}H_0^{(1)}\left(z\right)
\end{equation}
with
$$z=\frac{is}{a\sqrt{1-M^2}}\sqrt{\frac{(x-\xi)^2}{1-M^2}+y^2}$$
is well defined in $\R^2\setminus[-1,1]$ and represents a solution of equation
(\ref{RedWave}). In addition, formula
$$\frac{\partial}{\partial \eta}H_0^{(1)}\left(z\right)
=\frac{dH_0^{(1)}}{dz}(z)\cdot\frac{\partial z}{\partial \eta}
=\frac{dH_0^{(1)}}{dz}(z)\frac{-is(y-\eta)}{a\sqrt{1-M^2}}
\cdot\left(\sqrt{\frac{(x-\xi)^2}{1-M^2}+(y-\eta)^2}\right)^{-1}$$
shows that
\begin{equation}\label{DoubletZero}
\Psi_{\xi,0}(x,y)=0\hspace{0.05in}\mbox{for}\ \left\{|x|>1, y=0\right\}.
\end{equation}
\indent
Therefore, by Lemma~\ref{ReducedWave} function
$$\psi_{\xi,0}(x,y,t)=\Psi_{\xi,0}(x,y)e^{s(t+cx)},$$
with $c$ defined in formula (\ref{cFormula}), is a solution of equation (\ref{LinearEquation})
for arbitrary $s$. Following the spirit of terminology in \cite{BAH} we call this solution
a {\it $(\xi,s)$-doublet} and construct a {\it general doublet-solution}
of equation (\ref{LinearEquation}) by the formula
\begin{equation}\label{psiDoublet}
\psi(x,y,t)=\frac{1}{2\pi i}
\int_{\sigma-i\infty}^{\sigma+i\infty}dse^{s(t+cx)}
\int_{-1}^1d\xi p(\xi,s)\left(\frac{\partial H_0^{(1)}}{\partial \eta}(z)\Bigg|_{\eta=0}\right)
\end{equation}
with arbitrary function $p(\xi,s)$, satisfying condition
$$\left\|p(\xi,\sigma+i\nu)\right\|_{L^1_{\nu}}<K$$
uniformly with respect to $\xi$ and $\sigma\in[\sigma_1,\sigma_2]$.

\section{Generalized Possio integral equation.}\label{PossioEquation}

\indent
In this section we construct the velocity potential corresponding to a doublet-solution
in (\ref{psiDoublet}) and represent the flow tangency condition (\ref{FlowTangency})
as an integral equation with respect to $p(\xi,s)$, which we call {\it generalized Possio integral
equation}.
Since condition (\ref{DoubletZero}) is satisfied by any $\Psi_{\xi,0}$, condition (\ref{Kutta})
is automatically satisfied by arbitrary doublet-solution of the form (\ref{psiDoublet}).
Therefore, in order to solve the boundary value problem under consideration we only need
the corresponding velocity potential to satisfy condition (\ref{FlowTangency}).\\
\indent
We start the construction of the velocity potential with the construction of such a
potential for a $(\xi,s)$-doublet. Considering the definition of $\psi$
$$\psi=\frac{\partial\phi}{\partial t}+U\frac{\partial\phi}{\partial x},$$
as a differential equation with respect to $\phi$ and assuming that
$$\phi_{\xi,0}(x,y,t)=\Phi_{\xi,0}(x,y)e^{s(t+cx)},$$
we obtain an equation for $\Phi_{\xi,0}$
$$\frac{\partial\phi_{\xi,0}}{\partial t}+U\frac{\partial\phi_{\xi,0}}{\partial x}
=s\Phi e^{s(t+cx)}+Usc\Phi e^{s(t+cx)}+U\frac{\partial\Phi}{\partial x}e^{s(t+cx)}
=\frac{\partial}{\partial \eta}H_0^{(1)}(z)\Bigg|_{\eta=0}e^{s(t+cx)}$$
or
$$U\frac{\partial\Phi}{\partial x}+s\left(1+cU\right)\Phi
=\frac{\partial}{\partial \eta}H_0^{(1)}(z)\Bigg|_{\eta=0}.$$
\indent
Solving this equation we obtain
\begin{equation}\label{phiDoubletFormula}
\Phi_{\xi,0}(x,y)=\frac{e^{\lambda x}}{U}\int_{-\infty}^x e^{-\lambda u}
\left(\frac{\partial}{\partial \eta}H_0^{(1)}(\zeta)\Bigg|_{\eta=0}\right)du
\end{equation}
with
\begin{equation}\label{lambdaFormula}
\lambda=\lambda(s)=-\frac{s\left(1+cU\right)}{U},
\end{equation}
and
$$\zeta=\frac{is}{a\sqrt{1-M^2}}\sqrt{\frac{(u-\xi)^2}{1-M^2}+(y-\eta)^2}.$$
\indent
Finally, from (\ref{phiDoubletFormula}) using equality
$$\frac{\partial H_0^{(1)}}{\partial \eta}(z)\Bigg|_{\eta=0}
=\frac{dH_0^{(1)}}{dz}\cdot\frac{\partial z}{\partial\eta}\Bigg|_{\eta=0}
=-\frac{\partial H_0^{(1)}}{\partial y}(z)\Bigg|_{\eta=0}.$$
we obtain an expression for the velocity potential
corresponding to arbitrary doublet-solution
\begin{equation}\label{phiFormula}
\phi(x,y,t)=-\frac{1}{2\pi i}
\int_{\sigma-i\infty}^{\sigma+i\infty}dse^{s(t+cx)}\frac{e^{\lambda(s)x}}{U}
\int_{-1}^1p(\xi,s)d\xi\int_{-\infty}^x e^{-\lambda(s)u}
\left(\frac{\partial H_0^{(1)}}{\partial y}(\zeta)\Bigg|_{\eta=0}\right)du.
\end{equation}
\indent
Using formula (\ref{phiFormula}) we rewrite the flow tangency condition
$$\lim_{y\to 0^{+}}\frac{\partial\phi}{\partial y}(x,y,t)=w(x,t)\hspace{0.05in}\mbox{for}
\hspace{0.05in}x\in[-1,1]$$
as
$$w(x,t)$$
$$=-\frac{1}{2\pi i}\lim_{y\to 0^{+}}
\int_{\sigma-i\infty}^{\sigma+i\infty}dse^{s(t+cx)}\frac{e^{\lambda(s)x}}{U}
\int_{-1}^1p(\xi,s)d\xi\int_{-\infty}^x e^{-\lambda(s)u}
\left(\frac{\partial^2 H_0^{(1)}}{\partial y^2}(\zeta)\Bigg|_{\eta=0}\right)du.$$
\indent
Representing then
$$w(x,t)=\frac{1}{2\pi i}\int_{\sigma-i\infty}^{\sigma+i\infty}e^{st}
{\widehat w}(x,s)ds,$$
we transform the equation above into
$$\frac{1}{2\pi i}\int_{\sigma-i\infty}^{\sigma+i\infty}dse^{st}\left\{{\widehat w}(x,s)
+\lim_{y\to 0^{+}}e^{scx}\frac{e^{\lambda(s)x}}{U}
\int_{-1}^1p(\xi,s)d\xi\int_{-\infty}^x e^{-\lambda(s)u}\right.$$
$$\left.\times\left(\frac{\partial^2 H_0^{(1)}}{\partial y^2}(\zeta)
\Bigg|_{\eta=0}\right)du\right\}=0$$
and further into
$${\widehat w}(x,s)e^{-scx}=-\lim_{y\to 0^{+}}
\frac{e^{\lambda(s)x}}{U}\int_{-1}^1d\xi p(\xi,s)\int_{-\infty}^x due^{-\lambda(s)u}
\left(\frac{\partial^2 H_0^{(1)}}{\partial y^2}(\zeta)\Bigg|_{\eta=0}\right).$$
\indent
To transform equation above into a singular integral equation we follow two
additional steps from the Possio's scheme (\cite{BAH}). On the first step we use the fact
that $H_0^{(1)}(z)$ satisfies equation (\ref{RedWave}), and therefore
$$\frac{\partial^2 H_0^{(1)}}{\partial y^2}(z)
=-(1-M^2)\frac{\partial^2 H_0^{(1)}}{\partial x^2}(z)+\frac{s^2}{a^2(1-M^2)}H_0^{(1)}(z).$$
\indent
Then we obtain equation
$${\widehat w}(x,s)e^{-scx}=\lim_{y\to 0^{+}}
\frac{e^{\lambda(s)x}}{U}\int_{-1}^1d\xi p(\xi,s)$$
$$\times\int_{-\infty}^x due^{-\lambda(s)u}
\left((1-M^2)\frac{\partial^2 H_0^{(1)}}{\partial u^2}(\zeta)
-\frac{s^2}{a^2(1-M^2)}H_0^{(1)}(\zeta)\right)\Bigg|_{\eta=0}.$$
\indent
On the second step we transform the second integral in the equation above using integration
by parts in the formula
\begin{equation}\label{Hbyparts}
\int_{-\infty}^x e^{-\lambda u}
\frac{\partial^2 H_0^{(1)}}{\partial u^2}(\zeta)du
=e^{-\lambda u}\cdot\frac{\partial H_0^{(1)}}{\partial u}(\zeta)\Bigg|_{-\infty}^x
+\lambda\int_{-\infty}^x e^{-\lambda u}\frac{\partial H_0^{(1)}}{\partial u}(\zeta)du
\end{equation}
$$=e^{-\lambda x}\cdot\frac{\partial H_0^{(1)}}{\partial x}(z)
+\lambda e^{-\lambda u}\cdot H_0^{(1)}(\zeta)\Bigg|_{-\infty}^x
+\lambda^2\int_{-\infty}^x e^{-\lambda u}H_0^{(1)}(\zeta)du$$
$$=e^{-\lambda x}\cdot\frac{\partial H_0^{(1)}}{\partial x}(z)
+\lambda e^{-\lambda x}\cdot H_0^{(1)}(z)
+\lambda^2\int_{-\infty}^x e^{-\lambda u}H_0^{(1)}(\zeta)du.$$
\indent
Using formula (\ref{Hbyparts}) in the equation above we obtain the generalized Possio
integral equation
\begin{equation}\label{Possio}
{\widehat w}(x,s)e^{-scx}=\int_{-1}^1d\xi p(\xi,s)
\lim_{y\to 0^{+}}\left[\frac{(1-M^2)}{U}\cdot\frac{\partial H_0^{(1)}}{\partial x}(z)
+\frac{\lambda(s)(1-M^2)}{U}H_0^{(1)}(z)\right.
\end{equation}
$$\left.+\frac{e^{\lambda(s)x}}{U}\left(\lambda^2(s)(1-M^2)-\frac{s^2}{a^2(1-M^2)}\right)
\int_{-\infty}^x e^{-\lambda(s) u}H_0^{(1)}(\zeta)du\right]\Bigg|_{\eta=0}.$$
\indent
This equation reduces to the Possio integral equation (\cite{BAH}) if we fix $s$ and consider
oscillations with fixed frequency. In this case function $w(x,t)$ reduces to
$w(x)e^{ikt}$ and the sought solution $p(\xi,s)$ reduces to $p(\xi)$.

\section{Finite Hilbert transform.}\label{Finite}

\indent
Before analyzing the solvability of the integral equation (\ref{Possio}) we represent
this equation as an integral equation with the principal term being a multiple
of the finite Hilbert transform \cite{Tr}. In order to obtain such a representation
we use Taylor (respectively Laurent) series of Bessel functions $Y_0$ and $Y_1$
(\cite{EMOT}, v.II, 7.2.4, (32)).\\

\begin{proposition}\label{Representation}\ Kernel of the generalized Possio equation
(\ref{Possio}) admits the following representation
\begin{equation}\label{RRepresentation}
\lim_{y\to 0^{+}}\left[\frac{(1-M^2)}{U}\cdot\frac{\partial H_0^{(1)}}{\partial x}(z)
+\frac{\lambda(s)(1-M^2)}{U}H_0^{(1)}(z)\right.
\end{equation}
$$\left.+\frac{e^{\lambda(s)x}}{U}\left(\lambda^2(s)(1-M^2)-\frac{s^2}{a^2(1-M^2)}\right)
\int_{-\infty}^x e^{-\lambda(s) u}H_0^{(1)}(\zeta)du\right]\Bigg|_{\eta=0}$$
$$=-\frac{2i(1-M^2)^{3/2}}{\pi U}\frac{1}{(x-\xi)}+K(x,\xi,s),$$
where
\begin{equation}\label{KRepresentation}
K(x,\xi,s)=A(x,\xi,s)\log{|x-\xi|}+B(x,\xi,s),
\end{equation}
with functions $A$ and $B$ bounded with respect to $x,\xi$, analytically depending on $s$
and satisfying
estimate
\begin{equation}\label{ABEstimate}
\sup_{x,\xi\in[-1,1]}\left\{\left|A_1(x,\xi,s)\right|,\
\left|B_1(x,\xi,s)\right|\right\}<C|s|^2.
\end{equation}
\end{proposition}
\indent
{\bf Proof.}\ We consider the first term of the kernel of equation (\ref{Possio}),
for which using equality
$$\frac{d}{dz}H_0^{(1)}(z)=-H_1^{(1)}(z)$$
we obtain
\begin{equation}\label{Representation1}
\frac{1-M^2}{U}\lim_{y\to 0^{+}}
\left(\frac{\partial H_0^{(1)}(z)}{\partial x}\right)\Bigg|_{\eta=0}
=\frac{1-M^2}{U}\lim_{y\to 0^{+}}\frac{dH_0^{(1)}}{dz}(z)
\cdot\frac{\partial z}{\partial\eta}\Bigg|_{\eta=0}
\end{equation}
$$=-\frac{1-M^2}{U}H_1^{(1)}\left(\frac{is|x-\xi|}{a(1-M^2)}\right)
\cdot\frac{is(x-\xi)}{a|x-\xi|\sqrt{1-M^2}}.$$
$$=-\frac{1-M^2}{U}\left[J_1\left(\frac{is|x-\xi|}{a(1-M^2)}\right)
+iY_1\left(\frac{is|x-\xi|}{a(1-M^2)}
\right)\right]\cdot\frac{is(x-\xi)}{a|x-\xi|\sqrt{1-M^2}},$$
where $J_1$ and $Y_1$ are Bessel functions.\\
\indent
For the case $\left|s(x-\xi)\right|<C$ with $\mbox{Re}s\in [\sigma_1,\sigma_2]$
we use the Laurent series of $Y_1$, analyticity of $J_1$ and boundedness of
$$\frac{is(x-\xi)}{a|x-\xi|\sqrt{1-M^2}}$$
and obtain the following formula
\begin{equation}\label{FirstTerm1}
\frac{1-M^2}{U}\lim_{y\to 0^{+}}
\left(\frac{\partial H_0^{(1)}(z)}{\partial x}\right)\Bigg|_{\eta=0}
\end{equation}
$$=-\frac{2i(1-M^2)^{3/2}}{\pi U}\frac{1}{(x-\xi)}+A_1(x,\xi,s)\log{|x-\xi|}
+B_1(x,\xi,s)$$
with functions $A_1$ and $B_1$ bounded with respect to $x,\xi$ for fixed $s$,
analytically depending on $s$, and satisfying
\begin{equation}\label{ABEstimate1}
\sup_{x,\xi\in[-1,1]}\left\{\left|A_1(x,\xi,s)\right|,\
\left|B_1(x,\xi,s)\right|\right\}<C|s|
\end{equation}
for some $C>0$.\\
\indent
For the case $\left|s(x-\xi)\right|>C$ with $\mbox{Re}s\in [\sigma_1,\sigma_2]$
we again use formula (\ref{Representation1}) and asymptotic expansions of Hankel functions
for large $|z|$ (\cite{EMOT}, v.II, 7.13). Then we obtain representation
\begin{equation}\label{FirstTerm2}
\frac{1-M^2}{U}\lim_{y\to 0^{+}}
\left(\frac{\partial H_0^{(1)}(z)}{\partial x}\right)\Bigg|_{\eta=0}
=A_2(x,\xi,s)\log{|x-\xi|}+B_2(x,\xi,s)
\end{equation}
with functions $A_2$ and $B_2$ bounded with respect to $x,\xi$ for fixed $s$,
analytically depending on $s$, and satisfying estimate (\ref{ABEstimate1}).\\
\indent
For the rest of the kernel of equation (\ref{Possio}) in the case
$\left|s(x-\xi)\right|<C$ with $\mbox{Re}s\in [\sigma_1,\sigma_2]$ we use the
Taylor series of $Y_0$ and obtain
\begin{equation}\label{SecondTerm}
\begin{array}{ll}
{\dis \lim_{y\to 0^{+}}\left[\frac{\lambda(s)(1-M^2)}{U}H_0^{(1)}(z) \right.}
\vspace{0.1in}\\
{\dis \left.+\frac{e^{\lambda(s)x}}{U}\left(\lambda^2(s)(1-M^2)
-\frac{s^2}{a^2(1-M^2)}\right)
\int_{-\infty}^x e^{-\lambda(s)u}H_0^{(1)}(\zeta)du\right]\Bigg|_{\eta=0} }
\end{array}
\end{equation}
$$=A_3(x,\xi,s)\log{|x-\xi|}+B_3(x,\xi,s)$$
with functions $A_3(x,\xi,s)$ and $B_3(x,\xi,s)$
analytically depending on $s$ and satisfying estimate (\ref{ABEstimate}).\\
\indent
In the case$\left|s(x-\xi)\right|>C$ with $\mbox{Re}s\in [\sigma_1,\sigma_2]$ we
again use asymptotic expansions of Hankel functions for large $|z|$. Then we obtain
representation, similar to (\ref{SecondTerm}) with functions satisfying
estimate (\ref{ABEstimate}).\\
\indent
Combining formulas above we obtain representation (\ref{RRepresentation})
of the kernel of the generalized Possio equation with functions $A(x,\xi,s)$
and $B(x,\xi,s)$ satisfying estimate (\ref{ABEstimate}).\qed\\

\section{Solvability of the generalized Possio equation.}\label{Solvability}

\indent
Using representation (\ref{RRepresentation}) we consider operators
$${\cal K}_s[f](x)=\int_{-1}^1f(\xi)K(x,\xi,s)d\xi$$
and
$${\cal R}_s[f](x)={\cal T}[f](x)-\frac{iU}{2(1-M^2)^{3/2}}{\cal K}_s[f](x),$$
where ${\cal T}$ is the finite Hilbert transform \cite{Tr}:
$${\cal T}[f](x)=\frac{1}{\pi}\int_{-1}^1\frac{f(\xi)}{\xi-x}d\xi.$$
\indent
Then we can rewrite equation (\ref{Possio}) in the following form
\begin{equation}\label{REquation}
{\cal R}_s[p](x)=-\frac{iU}{2(1-M^2)^{3/2}}{\widehat w}(x,s)e^{-scx},
\end{equation}
with unknown function $p(\xi,s)$.\\
\indent
An important role in the analysis of solvability
of equation (\ref{Possio}) plays operator described in the following proposition (\cite{So}, \cite{Tr}).

\begin{proposition}\label{TInverse}\ Operator ${\cal T}^{-1}$ defined by the formula
\begin{equation}\label{T-1Operator}
{\cal T}^{-1}[g](x)=-\frac{1}{\pi}\int_{-1}^{1}\sqrt{\frac{1-y^2}{1-x^2}}
\frac{g(y)}{y-x}dy,
\end{equation}
is a bounded linear operator from $L^{\frac{4}{3}+}\left[-1,1\right]$ into $L^p\left[-1,1\right]$
for any $p<\frac{4}{3}$, satisfying equality
$${\cal T}\circ{\cal T}^{-1}[f]=f.$$
\end{proposition}
\qed\\
\indent
Using operator ${\cal T}^{-1}$ we reduce solution of equation (\ref{REquation}), and therefore
of equation (\ref{Possio}), to the solution of equation
\begin{equation}\label{GEquation}
{\cal G}_s[r](x)=-\frac{iU}{2(1-M^2)^{3/2}}{\widehat w}(x,s)e^{-scx},
\end{equation}
where
\begin{equation}\label{GOperator}
{\cal G}_s={\cal R}_s\circ {\cal T}^{-1}={\cal I}+{\cal N}_s,
\end{equation}
with ${\cal I}$ - the identity operator and
$${\cal N}_s=-\frac{iU}{2(1-M^2)^{3/2}}{\cal K}_s\circ{\cal T}^{-1}.$$
\indent
The advantage of equation (\ref{GEquation}) over equation (\ref{REquation}) becomes clear
from the proposition below, in which we prove the Fredholm property of the family
of operators ${\cal G}_s$. This proposition is the key new ingredient in the analysis of
solvability of the Possio integral equation and is inspired by the Proposition 5.1 from \cite{P}.\\

\begin{proposition}\label{GFredholm}\ For any fixed $s\in\C$ operator
${\cal N}_s$ is compact on $L^{2}[-1,1]$, and therefore operator ${\cal G}_s$ defined
in (\ref{GOperator}) is a Fredholm operator on $L^{2}[-1,1]$. In addition, kernel
$N(x,y,\lambda)$ of the operator ${\cal N}_s$ admits estimate
\begin{equation}\label{NEstimate}
\int_{\R^2}|N(x,y,s)|^2dxdy<C|s|^4
\end{equation}
with constant $C$ independent of $s$.
\end{proposition}
\indent
{\bf Proof.}\ Using formula (\ref{T-1Operator}) for ${\cal T}^{-1}$ we obtain
$${\cal N}_s[f](x)=-\frac{iU}{2(1-M^2)^{3/2}}{\cal K}_s\left[-\frac{1}{\pi}
\int_{-1}^{1}\sqrt{\frac{1-y^2}{1-u^2}}\frac{f(y)}{y-u}dy\right]$$
$$=\frac{iU}{2(1-M^2)^{3/2}}\int_{-1}^1duK(x,u,s)
\int_{-1}^{1}\sqrt{\frac{1-y^2}{1-u^2}}\frac{f(y)}{y-u}dy$$
$$=\int_{-1}^{1}N(x,y,s)f(y)dy,$$
where
$$N(x,y,s)=\frac{iU}{2(1-M^2)^{3/2}}\int_{-1}^{1}du\frac{K(x,u,s)}{y-u}\sqrt{\frac{1-y^2}{1-u^2}}.$$
\indent
To prove compactness of the operator ${\cal N}_s$ we use representation
$$N(x,y,s)=\frac{iU}{2(1-M^2)^{3/2}}\left[N_1(x,y,s)+N_2(x,y,s)\right],$$
with
$$N_1(x,y,s)=\int_{-1}^{1}K(x,u,s)\frac{du}{y-u},$$
and
$$N_2(x,y,s)=\int_{-1}^{1}K(x,u,s)
\frac{\left(\sqrt{1-y^2}-\sqrt{1-u^2}\right)}{\sqrt{1-u^2}}\frac{du}{y-u}$$
$$=-\int_{-1}^{1}K(x,u,s)
\frac{\left(y+u\right)du}{\left(\sqrt{1-y^2}+\sqrt{1-u^2}\right)\sqrt{1-u^2}},$$
and prove Hilbert-Schmidt property (cf.\cite{L})
of kernels $N_1(x,y,s)$ and $N_2(x,y,s)$.\\
\indent
For $N_1$ we notice that according to estimates (\ref{ABEstimate}) for fixed $x\in[-1,1]$
$$\int_{-1}^{1}K(x,u,s)\frac{du}{y-u}$$
is a multiple of the Hilbert transform of an $L^2[-1,1]$ - function
$K(x,u,s)$ satisfying
$$\|K(x,u,s)\|_{L^2_u}<C|s|^2$$
with constant $C$ independent of $x$. Therefore,
$$\int_{-1}^{1}dx\int_{-1}^{1}dy\left|N_1(x,y,s)\right|^2
=\int_{-1}^{1}dx\int_{-1}^{1}dy\left|\int_{-1}^{1}K(x,u,s)\frac{du}{y-u}\right|^2$$
$$<C\int_{-1}^{1}dx\left\|K(x,u,s)\right\|_{L^2[-1,1]_u}^2<C|s|^4.$$
\indent
For $N_2(x,y,s)$ we have
\begin{equation}\label{N2Integral}
\begin{array}{lll}
{\dis \int_{-1}^{1}dx\int_{-1}^{1}dy\left|N_2(x,y,s)\right|^2 }\vspace{0.1in}\\
{\dis =\int_{-1}^{1}dx\int_{-1}^{1}dy\left|\int_{-1}^{1}K(x,u,s)
\frac{\left(y+u\right)du}{\left(\sqrt{1-y^2}+\sqrt{1-u^2}\right)\sqrt{1-u^2}}\right|^2 }
\vspace{0.1in}\\
{\dis <C |s|^4\int_{-1}^{1}dx\int_{-1}^{1}dy\left|\int_{-1}^{1}
\frac{\log\left|x-u\right|du}{\left(\sqrt{1-y^2}+\sqrt{1-u^2}\right)\sqrt{1-u^2}}\right|^2, }
\end{array}
\end{equation}
\indent
where in the last inequality we used representation (\ref{KRepresentation}) of function $K$.\\
\indent
To estimate the last integral in (\ref{N2Integral}) we define
$$S_{x,y}=\left\{u:|x-u|\geq \frac{1}{2}|x-1|\cdot\sqrt{1-y^2}\right\}$$
and estimate separately integrals over $S_{x,y}$ and $[-1,1]\setminus S_{x,y}$.
For the integral over $S_{x,y}$, changing variable $u=\cos{\theta}$, we have
$$\int_{-1}^{1}dx\int_{-1}^{1}dy\left|\int_{-1}^{1}
\frac{\log\left|x-u\right|du}{\left(\sqrt{1-y^2}+\sqrt{1-u^2}\right)\sqrt{1-u^2}}\right|^2$$
$$<C\int_{-1}^{1}dx\int_{-1}^{1}dy\left|\int_{0}^{\pi}
\frac{\left(\log|1-x|+\log{\sqrt{1-y^2}}\right)d\theta}
{\left(\sqrt{1-y^2}+|\sin{\theta}|\right)}\right|^2$$
$$<C\int_{-1}^{1}dx\int_{-1}^{1}dy\left(\log|1-x|+\log{\sqrt{1-y^2}}\right)^2
\left(\log{\sqrt{1-y^2}}\right)^2<C.$$
\indent
For $u\in [0,1]\setminus S_{x,y}$ we have
$$|1-u|\geq|1-x|-|x-u|\geq|1-x|-\frac{1}{2}|1-x|\cdot\sqrt{1-y^2}\geq\frac{1}{2}|1-x|.$$
Therefore, for the integral over $[0,1]\setminus S_{x,y}$ we obtain
$$\int_{-1}^{1}dx\int_{-1}^{1}dy\left|\int_{[0,1]\setminus S_{x,y}}
\frac{\log\left|x-u\right|du}{\left(\sqrt{1-y^2}+\sqrt{1-u^2}\right)\sqrt{1-u^2}}
\right|^2$$
$$<C\int_{-1}^{1}dx\int_{-1}^{1}dy\left|\frac{1}{\sqrt{1-y^2}\sqrt{|1-x|}}
\int_{x-\frac{1}{2}|x-1|\cdot\sqrt{1-y^2}}^{x+\frac{1}{2}|x-1|\cdot\sqrt{1-y^2}}
\log\left|x-u\right|du\right|^2$$
$$<C\int_{-1}^{1}dx\int_{-1}^{1}dy
\left|\frac{\left(|1-x|\cdot\sqrt{1-y^2}\right)
\left(\log{|1-x|}+\log{\sqrt{1-y^2}}\right)}{\sqrt{1-y^2}\sqrt{|1-x|}}\right|^2<C.$$
\indent
The same estimate holds for the integral over $[-1,0]\setminus S_{x,y}$.\\
\indent
Combining the estimates above we obtain estimate (\ref{NEstimate}).
\qed

\indent
To formulate a criterion of solvability of equation (\ref{Possio}) for a fixed $s\in\C$
we define a complex number $s_0$ to be a {\it characteristic value of the family of
operators ${\cal G}_s$} if operator ${\cal G}_{s_0}$ is not invertible in $L^2[-1,1]$.
Using Propositions~\ref{TInverse} and \ref{GFredholm}
we obtain the following criterion.\\

\begin{proposition}\label{Criterion}\ If $s_0$ is not a characteristic value of the
family of operators ${\cal G}_s$ and ${\widehat w}(x,s_0)e^{-s_0cx}\in L^2[-1,1]$,
then there exists a function ${\dis p(\xi,s_0)\in L^{\frac{4}{3}-}[-1,1] }$
satisfying equation (\ref{Possio}).

\end{proposition}
\indent
{\bf Proof.}\ We consider a solution $r$ of equation
$${\cal G}_{s_0}[r]={\cal R}_{s_0}\circ{\cal T}^{-1}[r]
=-\frac{iU}{2(1-M^2)^{3/2}}{\widehat w}(x,s_0)e^{-s_0cx}.$$
Then we have for $p(\xi,s_0)={\cal T}^{-1}[r]$
$${\cal R}_{s_0}[p]={\cal R}_{s_0}\circ{\cal T}^{-1}[r]
=-\frac{iU}{2(1-M^2)^{3/2}}{\widehat w}(x,s_0)e^{-s_0cx}.$$
\qed

\section{ The resolvent of operator ${\cal G}_s$.}\label{GOperatorResolvent}

\indent
In this section we address the question of solvability of integral equation (\ref{GEquation})
with varying $s$. We construct the resolvent of the operator ${\cal G}_s$
and show that the resolvent is also a Fredholm operator analytically depending on
$s \in \left\{\mbox{Re}s>\sigma_1\right\}$.\\
\indent
Let ${\cal P}:L^2(\R)\to L^2(\R)$ be an integral operator with kernel $P(x,y)$ satisfying
Hilbert-Schmidt condition. Following \cite{C}, we consider for operator ${\cal P}$
Hilbert's modification of the original Fredholm's determinants:
$$D_{P,m}\left(t_1,\dots,t_m\right)
=\left|
\begin{array}{cccc}
0&P(t_1,t_2)&\cdots&P(t_1,t_m)\\
P(t_2,t_1)&0&\cdots&P(t_2,t_m)\\
\vdots&&&\vdots\\
P(t_m,t_1)&\cdots&P(t_m,t_{m-1})&0
\end{array}
\right|,$$
\begin{equation}\label{D_P}
D_{P}=1+\sum_{m=1}^{\infty}\delta_m
=1+\sum_{m=1}^{\infty}\frac{1}{m!}\int_{\R}\cdots\int_{\R}
D_{P,m}\left(t_1,\dots,t_m\right)dt_1\cdots dt_m,
\end{equation}
$$D_{P,m}\left(\begin{array}{c}
x\\
y
\end{array}t_1,\dots,t_m\right)
=\left|
\begin{array}{cccc}
P(x,y)&P(x,t_1)&\cdots&P(x,t_m)\\
P(t_1,y)&0&\cdots&P(t_1,t_m)\\
\vdots&&&\vdots\\
P(t_m,y)&\cdots&P(t_m,t_{m-1})&0
\end{array}
\right|,$$
and
\begin{equation}\label{D_Pfunction}
\begin{array}{ll}
{\dis D_{P}\left(\begin{array}{c}
x\\
y
\end{array}\right)
=P(x,y)+\sum_{m=1}^{\infty}\delta_m\left(\begin{array}{c}
x\\
y
\end{array}\right) }\vspace{0.1in}\\
{\dis =P(x,y)+\sum_{m=1}^{\infty}\frac{1}{m!}\int_{\R}\cdots\int_{\R}
D_{P,m}\left(
\begin{array}{c}
x\\
y
\end{array}t_1,\dots,t_m\right)dt_1\cdots dt_m. }
\end{array}
\end{equation}
\indent
We start with the following proposition, which summarizes the results from \cite{C}
(cf. also \cite{M}), that will be used in the construction of the resolvent of
${\cal G}_s$.

\begin{proposition}\label{Resolvent} (\cite{C})\ Let function
$P(x,y):\R^2\to \C$ satisfy Hilbert-Schmidt condition
$$\|P\|^2=\int_{\R^2}\left|P(x,y)\right|^2dxdy<\infty.$$
\indent
Then function
$D_{P}\left(\begin{array}{c}
x\\
y
\end{array}\right)\in L^2(\R^2)$ is well defined, and the following estimates hold:
\begin{equation}\label{DEstimate1}
\left|\delta_m\right|\leq\left(\frac{e}{m}\right)^{m/2}\|P\|^m,\
\left|D_{P}\right|\leq e^{\frac{\|P\|^2}{2}},
\end{equation}
\begin{equation}\label{DEstimate2}
\left|D_{P}\left(\begin{array}{c}
x\\
y
\end{array}\right)\right|
\leq e^{\frac{\|P\|^2}{2}}
\left(|P(x,y)|+\sqrt{e}\alpha(x)\beta(y)\right),
\end{equation}
where
$$\alpha^2(x)=\int_{\R}\left|P(x,t)\right|^2dt,\hspace{0.1in}
\beta^2(y)=\int_{\R}\left|P(t,y)\right|^2dt.$$
\indent
If $D_{P}\neq 0$
then kernel
\begin{equation}\label{HResolvent}
H(x,y)=\left[D_{P}\right]^{-1}
\cdot D_{P}\left(\begin{array}{c}
x\\
y
\end{array}\right)
\end{equation}
defines the resolvent of operator ${\cal I}-{\cal P}$, i.e. it satisfies the
following equations
\begin{equation}\label{ResolventEquations}
\begin{array}{ll}
{\dis H(x,y)+\int_{\R}P(x,t)\cdot H(t,y)dt=P(x,y), }
\vspace{0.1in}\\
{\dis H(x,y)+\int_{\R}P(t,y)\cdot H(x,t)dt=P(x,y), }
\end{array}
\end{equation}
and therefore operator ${\cal I}-{\cal H}$ is the inverse of
operator ${\cal I}+{\cal P}$.
\end{proposition}
\qed\\

\indent
In the next proposition we obtain necessary estimate for the resolvent ${\cal I}-{\cal H}_s$
with respect to $s$.

\begin{proposition}\label{GResolvent} (\cite{P})\ The set of characteristic values of the family
of operators ${\cal G}_s$ coincides with the set
$$E({\cal G})=\left\{s\in \C:\ \mbox{Re}s>\sigma_1,\ D_{N_s}=0\right\}$$
and consists of at most countably many isolated points.\\
\indent
For $s \notin E({\cal G})$ there exists an operator ${\cal H}_s$ with kernel
$H(x,y,s)$ satisfying the Hilbert-Schmidt condition and such that
operator ${\cal I}-{\cal H}_s$ is the inverse of operator ${\cal I}+{\cal N}_s={\cal G}_s$.\\
\indent
If function $D_{N}(s)=D_{N_s}$ has no zeros in a strip
$\left\{s=\sigma+i\nu:\ \sigma_1<\mbox{Re}s<\sigma_2\right\}$, then operator ${\cal H}_s$
admits estimate
\begin{equation}\label{HNormEstimate}
\left||{\cal H}_s\right\|
<\exp\left\{Ce^{|\nu|}\cdot(1+|\nu|)^4\right\}
\end{equation}
for some $C>0$, $s\in \left\{\sigma_1+\gamma <\mbox{Re}s <\sigma_2-\gamma\right\}$,
and arbitrary $\epsilon>0$.
\end{proposition}
\indent
{\bf Proof.}\ This proposition is a copy of Proposition 6.2 from \cite{P} except estimate
(\ref{HNormEstimate}). As in that proposition, analyticity of ${\cal H}_s$ and countability
of the set of characteristic values follow from Proposition~\ref{Resolvent} and
Theorem VI.14 from \cite{RS} applied to the family of operators ${\cal N}_s$.\\
\indent
Estimate (\ref{HNormEstimate}) is proved by applying the standard estimate
$$\left\|{\cal P}\right\|^2\leq \int_{\R^2}\left|P(x,y)\right|^2dxdy$$
for an integral operator ${\cal P}$ with kernel $P(x,y)$, estimate from the
Lemma~\ref{Estimate-ofD-1} below and the following estimate
\begin{equation}\label{D_NEstimate}
\left\|D_{N}\left(\left.\begin{array}{c}
x\\
y
\end{array}\right|s\right)\right\|
< \exp{\left\{C(1+|\nu|)^4\right\}}\cdot(1+|\nu|)^4\hspace{0.05in}\mbox{for}\ s=\sigma+i\nu.
\end{equation}
\indent
To obtain estimate (\ref{D_NEstimate}) we use estimate (\ref{DEstimate2}). Then we have
$$\left\|D_{N}\left(\left.\begin{array}{c}
x\\
y
\end{array}\right|s\right)\right\|^2< Ce^{\|{\cal N}_s\|^2}
\cdot\int_{\R^2}dxdy\left[|N_s(x,y)|^2+\alpha^2(x)\beta^2(y)\right]$$
$$=Ce^{\|{\cal N}_s\|^2}
\cdot\left[\int_{\R^2}|N_s(x,y)|^2dxdy+\int_{\R}\alpha^2(x)dx\cdot\int_{\R}\beta^2(y)dy\right]$$
$$=Ce^{\|{\cal N}_s\|^2}
\cdot\left[\int_{\R^2}|N_s(x,y)|^2dxdy+\int_{\R}dx\int_{\R}dt|N_s(x,t)|^2
\cdot\int_{\R}dy\int_{\R}dt|N_s(t,y)|^2\right]$$
$$<C\exp\left\{|s|^4\right\}\cdot\left(|s|^4+|s|^8\right)$$
where in the last inequality we used estimate (\ref{NEstimate}).\\
\indent
Finally, in order to use formula (\ref{HResolvent}) for the estimate of $\left\|{\cal H}_s\right\|$
we need an estimate of the function $\left|1/D_{N}(s)\right|$,
which is given in the lemma below.

\begin{lemma}\label{Estimate-ofD-1} (\cite{P})\ If function $D_{N}(s)=D_{N_s}$
has no zeros in the strip
$\left\{s:\ \sigma_1<\mbox{Re}s<\sigma_2\right\}$, then estimate
\begin{equation}\label{D-1Estimate}
\left|1/D_{N}(s)\right|<\exp\left\{Ce^{|\nu|}\cdot(1+|\nu|)^4\right\}
\end{equation}
holds for $s=\sigma+i\nu\in \left\{\sigma_1+\gamma <\mbox{Re}s <\sigma_2-\gamma\right\}$
with fixed $\gamma>0$ and arbitrary $\epsilon>0$.
\end{lemma}
\indent
{\bf Proof.}\ We consider a biholomorphic map
$$\Psi:\ \left\{s:\ \sigma_1<\mbox{Re}s<\sigma_2\right\}\to
\D(1)=\left\{z\in\C: |z|<1\right\},$$
defined by the formula
$$\Psi(s)=\frac{e^{i(s-\sigma_1)\frac{\pi}{\sigma_2-\sigma_1}}-i}
{e^{i(s-\sigma_1)\frac{\pi}{\sigma_2-\sigma_1}}+i}.$$
\indent
Denoting
$$w=u+iv=e^{i(s-\sigma_1)\frac{\pi}{\sigma_2-\sigma_1}},$$
we obtain for the circle $C(r)=\left\{z:\ |z|=r\right\}$
$$\Psi^{-1}\left(C(r)\right)=\left\{\sigma+i\nu:\
\left|e^{i(s-\sigma_1)\frac{\pi}{\sigma_2-\sigma_1}}-i\right|
=r\left|e^{i(s-\sigma_1)\frac{\pi}{\sigma_2-\sigma_1}}+i\right|\right\}$$
$$=\left\{u+iv:\ \left(u^2+v^2-2v+1\right)
=r^2\left(u^2+v^2+2v+1\right)\right\}$$
$$=\left\{u+iv:\ u^2+\left(v-\frac{1+r^2}{1-r^2}\right)^2
=\frac{4r^2}{(1-r^2)^2}\right\}.$$
\indent
Introducing coordinates
$$\rho=\mbox{Re}\frac{\pi(s-\sigma_1)}{\sigma_2-\sigma_1},\
\tau=\mbox{Im}\frac{\pi(s-\sigma_1)}{\sigma_2-\sigma_1},$$
such that
$$w=u+iv=e^{i(s-\sigma_1)\frac{\pi}{\sigma_2-\sigma_1}}
=e^{i\rho-\tau}=e^{-\tau}\left(\cos{\rho}+\sin{\rho}\right),$$
we can rewrite the last condition as a quadratic equation with respect
to $e^{-\tau}$ for fixed $\rho$
$$\left(e^{-\tau}-\sin{\rho}\frac{1+r^2}{1-r^2}\right)^2
+\cos^2{\rho}\left(\frac{1+r^2}{1-r^2}\right)^2-\frac{4r^2}{(1-r^2)^2}=0.$$
\indent
Solving equation above we obtain
$$e^{-\tau}=\sin{\rho}\frac{1+r^2}{1-r^2}
\pm\sqrt{\frac{4r^2}{(1-r^2)^2}-\cos^2{\rho}\left(\frac{1+r^2}{1-r^2}\right)^2}$$
with solutions existing for $\rho$ such that
$$|\cos{\rho}|\leq \frac{2r}{1-r^2}\frac{1-r^2}{1+r^2}=\frac{2r}{1+r^2}.$$
\indent
The maximal value for $e^{-\tau}$ is achieved at $\rho=\frac{\pi}{2}$ and it is
$$e^{-\tau}=\frac{1+r^2}{1-r^2}+\frac{2r}{1-r^2}
=\frac{1+r^2+2r}{1-r^2}=\frac{(1+r)^2}{1-r^2}=\frac{1+r}{1-r}.$$
Therefore the maximal value for $|s|$ is achieved at $\rho=\frac{\pi}{2}$, is
equal to $|s|=\log\left(\frac{1+r}{1-r}\right)$, and for $r=1-\delta$
we have the maximal value
\begin{equation}\label{s_max}
\max|s|=\log\left(\frac{1+r}{1-r}\right)=-\log{\delta}+\log{(2-\delta)}.
\end{equation}
\indent
Since function $D_{N}\left(s\right)$ has no zeros
in $\left\{s:\ \sigma_1<\mbox{Re}s<\sigma_2\right\}$ we can consider
analytic function $\log\left(D_{N}\left(s\right)\right)$
in this strip, and using estimates (\ref{DEstimate1}) and (\ref{NEstimate}),
and equality (\ref{s_max}), we obtain the following estimate for
$z=(1-\delta)e^{i\theta}$
$$\log\left|D_{N}\left(\Psi^{-1}(z)\right)\right|
\leq \frac{\left\|N_{\Psi^{-1}(z)}\right\|^2}{2}
\leq C\left|\Psi^{-1}(z)\right|^4
\leq C|\log{\delta}|^4.$$
\indent
Using then the Borel-Caratheodory inequality (\cite{Ti1}, \cite{Boa}) on disks with radii
$$1-2\delta=r<R=1-\delta,$$
we obtain
$$\left|\log\left(D_{N}\left(\Psi^{-1}(z)\right)\right)
\right|_{\{|z|=1-2\delta\}}$$
$$\leq\frac{2-4\delta}{\delta}\max_{|z|=R}
{\mbox{Re}\left\{\log\left(D_{N}\left(\Psi^{-1}(z)\right)\right)\right\}}
+\frac{1-\delta+1-2\delta}{\delta}|\log\left(D_{N}\left(\Psi^{-1}(0)\right)\right)|$$
$$< \frac{C}{\delta}\log^4{\delta},$$
or
$$-\frac{C}{\delta}\log^4{\delta}
<\mbox{Re}\left\{\log\left(D_{N}\left(\Psi^{-1}(z)\right)\right)
\right\}\Big|_{\{|z|=1-2\delta\}}
<\frac{C}{\delta}\log^4{\delta}.$$
\indent
From the last estimate we obtain an estimate for
the function $\left|1/D_{N}\left(\Psi^{-1}(z)\right)\right|$ in the disk
$\D(1-2\delta)$:
\begin{equation}\label{D-lowestimate}
\left|1/D_{N}\left(\Psi^{-1}(z)\right)\right|\Bigg|_{\{|z|\leq 1-2\delta\}}
< C\exp\left\{\frac{|\log{\delta}|^4}{\delta}\right\}
\end{equation}
for arbitrary $\epsilon>0$.\\
\indent
For a fixed $\rho\in\left(0,\pi\right)$ and arbitrary $\tau$ we have that
$\rho+i\tau\in\Psi^{-1}\left(\D(r)\right)$ with $r=1-2\delta$ if
$$e^{|\tau|}\leq \sin{\rho}\cdot\frac{1+r^2}{1-r^2}
+\sqrt{\frac{4r^2}{(1-r^2)^2}-\cos^2{\rho}\cdot\left(\frac{1+r^2}{1-r^2}\right)^2}$$
$$=\sin{\rho}\cdot\frac{2-4\delta+4\delta^2}{2\delta(2-2\delta)}
+\frac{\sqrt{4(1-2\delta)^2-\cos^2{\rho}\cdot(2-4\delta+4\delta^2)^2}}
{2\delta(2-2\delta)},$$
and therefore for any interval $\left[\gamma^{\prime},\pi-\gamma^{\prime}\right]$
there exist constants $C_1,\ C_2$ such that conditions
$$\rho\in \left[\gamma^{\prime},\pi-\gamma^{\prime}\right],\hspace{0.1in}
\frac{C_1}{\delta}<e^{|\tau|}<\frac{C_2}{\delta}$$
imply that $\rho+i\tau \in \Psi^{-1}\left(\D(1-2\delta)\right)$.\\
\indent
Using then estimate (\ref{D-lowestimate}) we obtain for $s$ with
$\mbox{Re}s \in \left[\sigma_1+\frac{\gamma^{\prime}(\sigma_2-\sigma_1)}{\pi},
\sigma_2-\frac{(\pi-\gamma^{\prime})(\sigma_2-\sigma_1)}{\pi}\right]$ estimate
$$\left|1/D_{N}\left(s\right)\right|
<\exp\left\{Ce^{|s|}\cdot(1+|s|)^4\right\}$$
for arbitrary $\epsilon>0$, which leads to estimate (\ref{D-1Estimate}).
\qed

\indent
Combining now estimate (\ref{D_NEstimate}) for ${\dis \left\|D_{N}\left(\left.\begin{array}{c}
x\\
y
\end{array}\right|s\right)\right\| }$ with estimate (\ref{D-1Estimate}) we obtain
estimate (\ref{HNormEstimate}).
\qed

\section{Proof of Theorem~\ref{Main}.}\label{Proof}

\indent
According to Lemmas~\ref{ReducedWave} and \ref{HankelSolution} any function defined
by formula (\ref{Solution}) satisfies equation (\ref{LinearEquation}). Also, as we pointed out earlier,
this function will automatically satisfy boundary condition (\ref{Kutta}). Therefore,
to prove Theorem~\ref{Main} we have to solve integral equation (\ref{Possio}) for the function
$p(\xi,s)$, so that the boundary condition (\ref{FlowTangency}) will be also satisfied.
As in Proposition~\ref{Criterion} in order to solve equation (\ref{Possio}) we first solve equation
(\ref{GEquation})
$${\cal G}_s[r_s]={\cal R}_{s}\circ{\cal T}^{-1}[r_s]
=-\frac{iU}{2(1-M^2)^{3/2}}{\widehat w}(x,s)e^{-scx}$$
and then define
$$p(\xi,s)={\cal T}^{-1}[r_s].$$
\indent
Using Proposition~\ref{Resolvent} we obtain that the solution of (\ref{GEquation})
is defined by the formula
$$r_s=\left({\cal I}-{\cal H}_s\right)
\left[-\frac{iU}{2(1-M^2)^{3/2}}{\widehat w}(x,s)e^{-scx}\right],$$
where ${\cal H}_s$ is an integral operator with kernel
$$H_s(x,y)=\left[D_{N_s}\right]^{-1}
\cdot D_{N}\left(\left.\begin{array}{c}
x\\
y
\end{array}\right|s\right)$$
admitting estimate (\ref{HNormEstimate}). Combining estimate (\ref{HNormEstimate})
with estimate (\ref{wCondition}) we obtain estimate
\begin{equation}\label{rEstimate}
\left\|r(\cdot,\sigma+i\nu)\right\|_{L^2[-1,1]}<C
\end{equation}
with constant $C$ independent of $s$. Then from Proposition~\ref{TInverse} we obtain
estimate (\ref{pEstimate}) for $p(\xi,s)={\cal T}^{-1}[r_s]$.
\qed

\end{document}